\documentclass[english,11pt]{amsart}
\usepackage{amsopn,amsthm,amsfonts,amsmath,amssymb,amscd}
\usepackage{babel}
\usepackage{amssymb}
\usepackage{enumerate}

\newtheorem{Theorem}{Theorem}[section]
\newtheorem{Lemma}[Theorem]{Lemma}
\newtheorem{Proposition}[Theorem]{Proposition}
\newtheorem{Corollary}[Theorem]{Corollary}

\newenvironment{Remark}{{\bf Remark.}}

\newcommand{\RR}{\mathbb{R}}
\newcommand{\ZZ}{\mathbb{Z}}

\newcommand{\TT}{\mathbb{T}}
\newcommand{\BB}{\mathcal{B}}

\newcommand{\mm}{\mathcal{M}}
\newcommand{\DEF}[1]{\emph{#1}}
\newcommand{\Gl}[2]{\mathop{GL}_{#1}({#2})}
\newcommand{\supp}{{\rm supp}}

\newcommand{\diam}[1]{{\rm diam}(#1)}

\begin{document}

\title{Diameters of commuting graphs of matrices over  semirings}
\author{David Dol\v zan, Damjana Kokol Bukov\v sek, Polona Oblak}
\date{\today}

\address{D.~Dol\v zan:~Department of Mathematics, Faculty of Mathematics
and Physics, University of Ljubljana, Jadranska 19, SI-1000 Ljubljana, Slovenia; e-mail: 
david.dolzan@fmf.uni-lj.si}
\address{D.~Kokol Bukov\v sek:~Department of Mathematics, Faculty of Mathematics
and Physics, University of Ljubljana, Jadranska 19, SI-1000 Ljubljana, Slovenia; e-mail: 
damjana.kokol@fmf.uni-lj.si}
\address{P.~Oblak: Faculty of Computer and Information Science,
Tr\v za\v ska 25, SI-1000 Ljubljana, Slovenia; e-mail: polona.oblak@fri.uni-lj.si}


\bigskip

\begin{abstract} 
We calculate the diameters of  commuting graphs of matrices over the binary Boolean semiring,
 the tropical semiring and an arbitrary nonentire commutative semiring. 
 We also find the lower bound for the diameter of the commuting graph of the semiring of matrices over an
 arbitrary  commutative entire antinegative  semiring.
\end{abstract}

\maketitle 
%\parindent=0cm

%-----------------------------------------------------
%-----------------------------------------------------
\section{Introduction}
%-----------------------------------------------------
%-----------------------------------------------------

\bigskip

%A semiring is an algebraic structure similar to a ring, but without the requirement that each element must have an additive inverse. 

{\bf Definition.}
A \emph{semiring} is a set $S$ equipped with binary operations $+$ and $\cdot$ such that $(S,+)$ is a commutative monoid with identity element 0, and $(S,\cdot)$ is a monoid with identity element 1. 
In addition, operations $+$ and $\cdot$ are connected by distributivity and 0 annihilates $S$. A semiring is 
\emph{commutative} if $ab=ba$ for all $a,b \in S$.

A semiring $S$ is called \DEF{antinegative}, if $a+b=0$ implies that $a=b=0$. Antinegative semirings
are also called \DEF{antirings} (or \DEF{zero-sum-free} semirings). A semiring is 
\DEF{entire} (or \DEF{zero-divisor-free})  if $ab=0$ implies that $a=0$ or $b=0$. 
A semiring is a \DEF{division} semiring if all nonzero elements have multiplicative inverses.

\bigskip

The simplest example of an antinegative semiring  is the \DEF{binary Boolean semiring}, the set $\{0,1\}$
with  $1+1=1\cdot 1=1$. 
We will denote the binary Boolean semiring by $\BB$.

Moreover, the set of nonnegative integers (or reals) with the usual operations of addition and multiplication,
is a commutative entire antinegative semiring, and all distributive lattices are antinegative semirings.

\smallskip

On the set $\RR \cup \{-\infty\}$, we define operations $a \oplus b = \max\{a,b\}$ and $a \odot b = a + b$,
where $-\infty + a = a + (-\infty) = -\infty$. It is easy to verify that ($\RR \cup \{-\infty\}, \oplus, \odot)$
is a semiring. It is denoted by $\TT$\ and called the \DEF{tropical semiring}, sometimes also the \emph{max-plus
semiring}.
Tropical semiring is a commutative entire antinegative division semiring and it is closely related to the 
\emph{max algebra}, i.e. the semiring of nonnegative reals $\RR_+$, where the addition is the same as
in the tropical semiring and the  multiplication is  the ordinary multiplication on reals.
Recently, commuting matrices in max algebra were studied in \cite{KatzSchnSerg10}.

Over past decades, the tropical semirings and other tropical structures  were widely investigated.
Let us mention only a few pioneering works \cite{AkiGauGut09, DevSanStur05, Simon94} that connect the theory of matrices over classical and
tropical worlds.

\bigskip

We denote by $\mm_n(S)$ the semiring of all $n \times n$ matrices over a semiring $S$ and by 
$\Gl{n}{S}$ the group of all invertible matrices in $\mm_n(S)$.
By $I_n \in \mm_n(S)$, we denote the identity matrix and by $0_n\in \mm_n(S)$ the zero matrix.
The matrix with the only nonzero entry 1 in the $i$-th row and the $j$-th column will be denoted by $E_{i,j}$.

If $X$ is a subset of $S$, let 
 $$C_S(X)=\{s\in S; \; sx=xs \text{ for all } x\in X\}$$ 
denote the \emph{centralizer} of $X$ in $S$. For $x \in S$ we also define $C_S(x)=C_S(\{x\})$.

\bigskip

For any subset $T$ of a semiring $S$, we denote by $\Gamma(T)$ the \DEF{commuting graph} of $T$. 
The vertex set $V(\Gamma(T))$ of
 $\Gamma(T)$ is the set of elements in $T \backslash C_S(T)$. 
 An unorderded pair of vertices $x-y$ is  an edge of $\Gamma(T)$ if $x \ne y$ and $xy=yx$.

 The sequence of edges $x_0 - x_1$, $ x_1 - x_2$, ..., $x_{k-1} - x_{k}$ is called \emph{a path of length $k$} and is denoted by $x_0 - x_1 - \ldots - x_k$. The \DEF{distance} between two vertices is the length of the shortest
path between them. The \DEF{diameter} of the graph is the longest distance between any two vertices of the 
graph $\Gamma$ and will be denoted by $\diam{\Gamma}$.

For example, the commuting graph of the set of all $2 \times 2$ nilpotent matrices over an entire antinegative
 semiring $S$, is a disconnected graph with two components (corresponding to the strictly upper, and strictly 
 lower triangular matrices), where both components are isomorphic to a complete graph on $\vert S \vert - 1$ 
 vertices.
 
 \bigskip
 
Recently, the commuting graphs of matrix rings  and semirings
(\cite{Abdollahi, AMRR06, AR06, DolKuzObl10, DolObl10, giudici, Mohammadian10}) 
and commuting graphs of various algebraic structures (\cite{IraJaf08, Segev01}) have been studied. 
They %commuting graphs 
give an illustrative way of describing centralizers of elements% in algebraic structures
.
It was proved in \cite[Cor.~7]{AMRR06} that the diameter of the commuting graph of the full matrix
ring over an algebraically closed field is equal to 4. For rings and fields that are not algebraically closed,
the commuting graph 
of the full matrix ring might not
be connected at all, or if it is connected, its diameter can be larger than 4. (See e.g.
\cite[Ex.~2.~15]{DolKuzObl10}, where it has been proven that $5\leq \diam{\Gamma(\mm_9(\ZZ_2))}<\infty$.)

In this paper, we continue with the investigation of the diameters of commuting graphs of full matrix semirings. 
In \cite[Prop.~10]{DolObl10} it has been shown that the diameter of the commuting graph of the full matrix 
semiring over % the binary 
%Boolean semiring 
$\BB$ is bounded between 3 and 4.  In the second section, we  prove that it is equal to 4. This implies that 
the diameter of  the commuting graph of $\mm_n(S)$ is at least 4 for every commutative entire
antinegative semiring $S$.
Using this, we prove in Theorem \ref{thm:mntropical} that  for the  tropical semiring $\TT$ the  diameter of 
$\Gamma(\mm_n(\TT))$ is equal to 4. 
In Section 4, we prove that  $\diam{\Gamma(\mm_n(S))} = 3$ for every nonentire commutative semiring $S$
  and all $n \geq 2$, thus generalizing  \cite[Thm.~1.1]{giudici}, where a similar result has been
  proven for $\mm_n(\ZZ_m)$, where $m$ is not a prime number.

\bigskip
\bigskip

%-----------------------------------------------------
%-----------------------------------------------------
\section{The commuting graph of the full matrix semiring over $\BB$}
%-----------------------------------------------------
%-----------------------------------------------------

\bigskip

In this section, we prove that the diameter of the commuting graph of the full $n\times n$ matrix semiring
 over Boolean semiring  $\BB$ is equal to 4 for all $n\geq 3$.

\bigskip

We start with a general lemma, which can be easily proved by a straightforward calculation.
Denote by $J_n$ the nilpotent $n \times n$ matrix $E_{1,2}+E_{2,3}+\ldots+E_{n-1,n}$.
 
 \bigskip
 
 \begin{Lemma}\label{lem:centrJ}
   If $S$ is a semiring and $n\geq 2$, then the centralizer of $J_n$  is equal to
   $
   C_{\mm_n(S)}(J_n)=\{a_0I_n+a_1J_n+a_2 J_n^2+\ldots+a_{n-1}J_n^{n-1}; \; a_i \in S\}
   $, and the centralizer of $J_n^T$  is equal to
   $C_{\mm_n(S)}(J_n^T)=\{b_0I_n+b_1J_n^T+b_2 \left(J_n^T\right)^2+\ldots+b_{n-1}\left(J_n^T\right)^{n-1}; \; b_i \in S\}
   $.\hfill$\blacksquare$
 \end{Lemma}

\bigskip

The following Theorem is the main result of this section.

\bigskip

\begin{Theorem}\label{thm:bbs}
 If $\BB$ is a binary Boolean semiring, then $\diam{\Gamma(\mm_2(\BB))}=\infty$ and 
 $$\diam{\Gamma(\mm_n(\BB))} =4$$ 
 for $n \geq 3$.
\end{Theorem}

\medskip

\begin{proof}
 It was proved in \cite[Prop.~10]{DolObl10} that $\diam{\Gamma(\mm_2(\BB))}=\infty$. 
 Let $n \geq 3$ and let $E \in \mm_n(\BB)$ be a matrix of all ones. 
 It was also proved in \cite[Prop.~10]{DolObl10} that the distance between any matrix $A \in \mm_n(\BB)$
 and matrix $E$ is at most 2. It follows that $\diam{\Gamma(\mm_n(\BB))} \le 4$ for $n \geq 3$.
 To prove the equality, we have to find two matrices with the distance between them at least 4. 
 
 Let $n = 3$ and let
 $$
 A = \left[\begin{matrix} 0 & 0 & 1 \\ 0 & 0 & 0 \\ 1 & 1 & 0 \end{matrix}\right] , \qquad 
 B = \left[\begin{matrix} 1 & 0 & 0 \\ 0 & 0 & 1 \\ 0 & 0 & 0 \end{matrix}\right] .
 $$
 A straight-forward calculation shows that the set of all matrices with distance 1 to matrix $A$ is 
 $$ 
 \mathbb{A} = \left\{ \left[\begin{matrix} 1 & 0 & 1 \\ 0 & 1 & 0 \\ 1 & 1 & 1 \end{matrix}\right] , 
        \left[\begin{matrix} 1 & 1 & 0 \\ 0 & 0 & 0 \\ 0 & 0 & 1 \end{matrix}\right] , 
        \left[\begin{matrix} 1 & 1 & 0 \\ 0 & 1 & 0 \\ 0 & 0 & 1 \end{matrix}\right] , 
        \left[\begin{matrix} 1 & 1 & 1 \\ 0 & 0 & 0 \\ 1 & 1 & 1 \end{matrix}\right] , 
        \left[\begin{matrix} 1 & 1 & 1 \\ 0 & 1 & 0 \\ 1 & 1 & 1 \end{matrix}\right] \right\} 
 $$
 and the set of all matrices with distance 1 to matrix $B$ is        
 $$ 
 \mathbb{B} = \left\{ \left[\begin{matrix} 1 & 0 & 0 \\ 0 & 0 & 0 \\ 0 & 0 & 0 \end{matrix}\right] , 
        \left[\begin{matrix} 1 & 0 & 0 \\ 0 & 1 & 1 \\ 0 & 0 & 1 \end{matrix}\right] , 
        \left[\begin{matrix} 0 & 0 & 0 \\ 0 & 0 & 1 \\ 0 & 0 & 0 \end{matrix}\right] , 
        \left[\begin{matrix} 0 & 0 & 0 \\ 0 & 1 & 0 \\ 0 & 0 & 1 \end{matrix}\right] , 
        \left[\begin{matrix} 0 & 0 & 0 \\ 0 & 1 & 1 \\ 0 & 0 & 1 \end{matrix}\right] \right\} .
 $$
 It is easy to check, that for any pair of matrices $C \in \mathbb{A}$ and $D \in \mathbb{B}$,
 $C$ and $D$ do not commute, so that the distance between $A$ and $B$ is at least 4
 and thus $\diam{\Gamma(\mm_3(\BB))}= 4$.
 
 Now, let $n \ge 4$. 
 Let $A, B \in \mm_n(\BB)$ be
 $$ 
 A = \left[\begin{matrix} 0 & \begin{matrix} 0 & \cdots & 0 & 1 \end{matrix}\\ \begin{matrix} 0 \\ 1 \\ \vdots \\ 1 \end{matrix} & J_{n-1}^T \end{matrix}\right] , \qquad 
 B = \left[\begin{matrix} 1 & 0  \\ 0 & J_{n-1} \end{matrix}\right] .
 $$
 Note that the centre of $\mm_n(\BB)$ consists only of $0_n$ and $I_n$. Suppose that the distance between 
 $A$ and $B$ is at most 3. Then there exist nonscalar matrices $C, D  \in \mm_n(\BB)$ such that 
 $A-C-D-B$ is a path in $\Gamma(\mm_n(\BB))$.
 
 Observe that $C$ is not a diagonal matrix, otherwise all of its diagonal entries 
 are equal, since $C$ commutes with $A$. So, $C + I$ is not in the centre and it commutes with
 $A$ and $D$. Since $C$ and $C + I$ have the same centralizer, we can therefore assume that all
 the  diagonal entries of $C$ are equal to $1$.
 
 Now, suppose that $D$ is diagonal. Since it commutes with $B$, it has the form
 $$ 
 D = \left[\begin{matrix} 1 & 0  \\ 0 & 0_{n-1} \end{matrix}\right] \qquad {\rm or} \qquad
 D = \left[\begin{matrix} 0 & 0  \\ 0 & I_{n-1} \end{matrix}\right] .
 $$
 In both cases, since $C$ and $D$ commute, $C$ has the form
 $$ 
 C = \left[\begin{matrix} 1 & 0  \\ 0 & C_1 \end{matrix}\right] .
 $$
 Since $C$ commutes with $A$, $C_1$ also commutes with $J_{n-1}^T$ and the last row of $C_1$ is equal to 
 $\left[\begin{matrix} 0 & \cdots & 0 & 1 \end{matrix}\right]$. So, $C_1 = I_{n-1}$ and $C$ is in the centre,
 which is a contradiction. Thus, $D$ is not diagonal and again we can assume that all diagonal entries of $D$ are 
 equal to $1$.
 
 Now, since $D$ and $B$  commute, we have 
 $$
 C = \left[\begin{matrix} 1 & c_{1,2} & \cdots & c_{1,n} \\ 
                          c_{2,1} & 1 & \cdots & c_{2,n} \\ 
                          \vdots & \vdots & \ddots & \vdots \\
                          c_{n,1} & c_{n,2} & \cdots & 1 \end{matrix}\right]
  \qquad {\rm and} \qquad                          D = \left[\begin{matrix} 1 & 0 &   0 & \cdots & 0       & 0 \\ 
                          0 & 1 & d_2 & \cdots & d_{n-2} & d_{n-1} \\ 
                          0 & 0 &   1 & \ddots &         & d_{n-2} \\
                          \vdots & \vdots &    & \ddots & \ddots & \vdots \\

                          0 & 0 &   0 & \cdots & 1        & d_2 \\
                          0 & 0 &   0 & \cdots & 0        & 1 \end{matrix}\right],
 $$
 where $d_i \ne 0$ for some $i \in \{2,3,\ldots,n-1\}$.
 
 Since $A$ and $C$ commute, we have for every $i = 2, ..., n-1$ that
 $ 0 = (AC)_{2,i} = (CA)_{2,i} = c_{2, i+1}$ and $0 = (AC)_{2,n} = (CA)_{2,n} = c_{2,1}$, so 
\begin{equation}\label{eq:0}
 c_{2,i} = 0 \quad {\rm for\ every} \quad i \ne 2.
\end{equation}
 Also, for every $i = 2, ..., n-1$ we have $(AC)_{1,i} = (CA)_{1,i}$ and thus
 \begin{equation}\label{eq:1}
 c_{n,i} = c_{1, i+1}.
 \end{equation} 
 Further, for every $i = 3, ..., n-1$ we have
 $ c_{1,i} = c_{1,i} + c_{2,i} = (AC)_{3,i} = (CA)_{3,i} = c_{3, i+1}$ and 
 $c_{1,n} = (AC)_{3,n} = (CA)_{3,n} = c_{3,1}$.
For  $i \ge j \geq 3$ and using 
 \begin{equation}\label{eq:4}
c_{1,i} + c_{j-1, i} = (AC)_{j,i} = (CA)_{j,i} = c_{j, i+1},
\end{equation}
 we prove by induction that 
 \begin{equation}\label{eq:2}
  c_{j, i+1} = c_{1,i} + c_{1, i-1} + ... + c_{1, i-j+3} \quad {\rm for\ every} \quad i \ge j \ge 3
 \end{equation} 
 and
 \begin{equation}\label{eq:3}
  c_{j,1} = c_{1,n} + c_{1, n-1} + ... + c_{1, n-j+3} \quad {\rm for\ every} \quad j \ge 3 .
   \end{equation} 
 Let $3 \le i < j \le n-1$. First $c_{j,3} = (CA)_{j,2} = (AC)_{j,2} = c_{1,2} + c_{j-1,2}$ and
 $c_{j,i} = (CA)_{j,i-1} = (AC)_{j,i-1} = c_{1, i-1} + c_{j-1, i-1}$ and by induction we have
 $$ c_{j,i} = c_{1,2} + ... + c_{1, i-1} +  c_{j-i+2,2} \quad {\rm for\ every} \quad 3 \le i < j \le n-1.$$ 
 
 Since $D$ is not diagonal, at least one of $d_2, ..., d_{n-1}$ is nonzero.
 Let $k$ be the greatest index, such that $d_k\ne 0$, so $d_k=1$. Note that the second row of $D$ is
 therefore  equal to
  $\left[\begin{matrix} 0 & 1 & d_2 & \cdots & d_{k-1} & 1 & 0 & \cdots & 0 \end{matrix}\right]$.
 
 Suppose first that $k = n-1$. Since $C$ and $D$ commute, we have 
 $0 = (CD)_{2,1} = (DC)_{2,1} = c_{2,1} + 
 d_2 c_{3,1} + ... + d_{n-1}c_{n,1} = c_{2,1} + 
 d_2 c_{3,1} + ... + d_{n-2}c_{n-1,1}+c_{n,1}$, thus $c_{n,1}=0$
 and by \eqref{eq:3}, we have $  c_{1,3} + ... + c_{1,n}=c_{n,1}=0$. It follows that
 $c_{1,i} = 0$ for all $i \ge 3$. Furthermore, $0 = c_{1,n} = (DC)_{1,n} = (CD)_{1,n} = c_{1,2}d_{n-1} = c_{1,2}$
 so  $c_{1,2} = 0$ and by \eqref{eq:2} it follows that $C$ is lower-triangular.
 
 Suppose now that $k < n-1$.
 Since $C$ and $D$ commute, we have
  $0 = (CD)_{2,1} = (DC)_{2,1} = c_{2,1} +   d_2 c_{3,1} + ... + d_{k}c_{k+1,1} = 
  c_{2,1} +   d_2 c_{3,1} + ... + d_{k-1}c_{k,1}+ c_{k+1,1}$, thus 
  $ c_{k+1,1} =0$ and by \eqref{eq:3} it follows that $ c_{1, n-k+2} + ... + c_{1,n}=0$. Therefore, 
 $c_{1,i} = 0$ for every $i \ge n-k+2$. Furthermore, 
 $0 = c_{1,n} = (DC)_{1,n} = (CD)_{1,n} = c_{1,n-k+1}d_{k} = c_{1,n-k+1}$.
 By induction, $c_{1,n-k-i+2} = c_{1,n-k-i+2}d_{k} = (CD)_{1,n-i+1} = (DC)_{1,n-i+1} = c_{1,n-i+1} = 0$
 for every $i= 1, ..., n-k$, so that $c_{1,j} = 0$ for every $j \ge 2$ and by \eqref{eq:2},
 $C$ is lower-triangular.
 
Now, we have $c_{j,i} = 0$ for every $i>j$ and by applying  \eqref{eq:4} it follows that 
\begin{equation}\label{eq:5}
 c_{j,i}=c_{j-1,i-1}=\ldots = c_{j-i+2,2}
\end{equation}
 for every $3 \le i < j \le n-1$. Since $C$ is lower-triangular, \eqref{eq:1} implies that $c_{n,i}=0$
 for $i\geq 2$ and $c_{n,1}=0$ by \eqref{eq:3}. Therefore we have for every $i = 2, ..., n-2$ that
 $ 0 = c_{n,i+1}=(CA)_{n,i} = (AC)_{n,i} = c_{n-1, i} $, which by \eqref{eq:5} implies that  $c_{n-i+1, 2}=0$,
 so again by \eqref{eq:5}, $c_{j,i}=0$ for all $3\leq i < j \leq n-1$. Together with \eqref{eq:0} we have that
  $C=I_n$, a contradiction. 
  
  So, we have proved that the distance between $A$ and $B$ cannot be less than 4 
  and therefore $\diam{\Gamma(\mm_n(\BB))}=4$.
\end{proof}

\bigskip

\begin{Corollary}\label{thm:ceas}
   If $S$ is a commutative entire antinegative semiring, then  $\diam{\Gamma(\mm_2(S))}=\infty$ and 
 $\diam{\Gamma(\mm_n(S))} \geq 4$, for $n \geq 3$.
\end{Corollary}
\medskip

\begin{proof}
  For a matrix $A \in \mm_n(S)$, let us denote by $\supp(A) \in \mm_n(\BB)$ the unique $(0,1)$-matrix 
  with the property
  $A_{i,j}\ne 0$ if and only if $\left(\supp(A)\right)_{i,j} \ne 0$ for all $1 \leq i,j \leq n$.
  
  Since $S$ is a commutative entire antinegative semiring, $AB=BA$ for $A, B \in \mm_n(S)$ implies that
  $\supp(A)\supp(B)=\supp(B)\supp(A)$ for $\supp(A), \supp(B) \in \mm_n(\BB)$. 
  Thus, it follows that $\diam{\Gamma(\mm_n(S))} \geq \diam{\Gamma(\mm_n(\BB))}$.
  Now, the statement follows by Theorem \ref{thm:bbs}.
\end{proof}

\bigskip
\bigskip

%-----------------------------------------------------
%-----------------------------------------------------
\section{The commuting graph of the full matrix semiring over $\TT$}
%-----------------------------------------------------
%-----------------------------------------------------

 \bigskip

The two operations $\oplus$ and $\odot$ in $\TT$ naturally induce the matrix addition and multiplication on 
the semiring $(\mm_n(\TT),\oplus,\odot)$, namely for $A=[a_{i,j}], B=[b_{i,j}] \in \mm_n(\TT)$ we have 
 \begin{align*}
   (A\oplus B)_{i,j}&=a_{i,j}\oplus b_{i,j}, \text{ and}\\
   (A \odot B)_{i,j}&=a_{i,1}\odot b_{1,j} \oplus a_{i,2}\odot b_{2,j}\oplus\ldots\oplus a_{i,n}\odot b_{n,j}.
 \end{align*}
 Let $I_{n}$ be the identity matrix in $\mm_n(\TT)$, i.e. matrix with zeros on the diagonal and $- \infty$ offdiagonal,
 and  let $E$ be the matrix with all entries equal to $0 \in \RR$. For $a\in \TT$ and $A=[a_{i,j}] \in \mm_n(\TT)$ we 
 also define the matrix $aA$ in the natural way, i.e. $(aA)_{i,j}=a\odot a_{ij}=a+a_{ij}$.

\bigskip

 It was proved in \cite[Cor.~11]{DolObl10} that  $\diam{\Gamma(\mm_n(\TT))} \geq 3$ for all $n \geq 3$.
Here, we will prove that $\diam{\Gamma(\mm_n(\TT))} = 4$ for all $n \geq 3$.

\bigskip
\begin{Lemma}\label{thm:0_n}
 The centralizer $C_{\mm_n(\TT)}(E)$ consists of exactly all matrices $A=[a_{i,j}]$, such that there 
 exists some element $a\in \TT$, with the property
   $$\max\limits_{j}\{a_{i,j}\}=a \text{ for all } i=1,2,\ldots,n \text{ and } \max\limits_{i}\{a_{i,j}\}=a 
   \text{ for all } j=1,2,\ldots,n.$$
\end{Lemma}

\medskip

\begin{proof}
 For the matrix $A=[a_{i,j}]$ we denote $a = \max_{i,j = 1, ..., n} \{a_{i,j}\}$.
 If $A$ commutes with $E$,  then 
 $\max\{a_{i,1}, a_{i,2}, ..., a_{i,n}\} = a_{i,1} \oplus a_{i,2} \oplus ... \oplus a_{i,n}
 = (A \odot E)_{i,1} = (E \odot A)_{i,1} = a_{1,1} \oplus a_{2,1} \oplus ... \oplus a_{n,1} =
  \max\{a_{1,1}, a_{2,1}, ..., a_{n,1}\}$
 for every $i = 1, ..., n$. So all rows of $A$ have the same maximum $a$, thus 
 $a$ appears in every row of $A$.
 Similarly, $a$ appears in every column of $A$. Conversely, if $a$ appears in every row and column of $A$, 
 then  clearly $E \odot A = a  E= A \odot E$. 
\end{proof}
 
\bigskip

\begin{Theorem}\label{thm:mntropical}
For the tropical semiring $\TT$, we have 
 $\diam{\Gamma(\mm_2(\TT))}=\infty$ and 
 $$\diam{\Gamma(\mm_n(\TT))} =4$$ for $n \geq 3$.
\end{Theorem}

\medskip

\begin{proof}
By Corollary \ref{thm:ceas} we have that $\diam{\Gamma(\mm_2(\TT))}=\infty$ and 
$\diam{\Gamma(\mm_n(\TT))} \geq 4$.

Suppose now $n \geq 3$. If $D=[d_{k,l}]$ is a diagonal matrix, where $d_{i,i}=d_{j,j}$ for some $i, j$, then let us define the $n \times n$ matrix $F=[f_{k,l}]$, such that 
 $$f_{k,l}=\begin{cases}
  0, & k=l \text{ or } (k,l)=(i,j)  \text{ or } (k,l)=(j,i) ,\\
  -\infty, & \text{otherwise} \, .
 \end{cases}$$
Now, $  (F \odot D)_{j,i}=(F \odot D)_{i,j}=d_{j,j}=d_{i,i}=(D \odot F)_{i,j} = (D \odot F)_{j,i}$. Furthermore,
$(F \odot D)_{k,k}=d_{k,k}=(D \odot F)_{k,k}$ for all $k$ and $(F \odot D)_{k,l}=-\infty=(D \odot F)_{k,l}$
for all $k \ne l$, and $\{k,l\} \ne \{i,j\}$. 
Thus $D-F-E$ is a path in $\Gamma(\mm_n(\TT))$
by Lemma \ref{thm:0_n}.

If $A=[a_{i,j}]$ is a nondiagonal matrix, let $a=\max\{a_{i,j}\}$, and then, $A-A\oplus a I_n-E$ is a path
in $\Gamma(\mm_n(\TT))$ by Lemma \ref{thm:0_n}.  

We have thus proved that $d(A,E) \leq 2$ for all matrices $A$ except for the diagonal matrices with all 
 of their diagonal entries distinct. 
We will now prove that $d(D,A) \leq 4$ for a diagonal matrix $D$
with all diagonal entries  distinct  and an arbitrary $A \in \mm_n(\TT)$. 

If $A$ is diagonal as well, it is clear that $d(D,A)=1$. So, suppose
$A=[a_{i,j}]$ is nondiagonal and let $\mu=\max\{a_{i,j}\}$ and $\varepsilon=\min\{a_{i,j}; \; a_{i,j} \ne -\infty\}$. 
Denote by $B=[b_{i,j}]$ the $n \times n$ matrix, defined by
 $$b_{i,j}=\begin{cases}
  \mu, & i=j,\\
  \varepsilon, & i \ne j\, .
 \end{cases}$$
Note that $\varepsilon + a_{i,j} \leq \varepsilon + \mu \leq a_{k,l} + \mu$ for all $i,j,k,l$ and thus for every
$1 \leq i,j\leq n$ we have $(A\odot B)_{i,j}=a_{i,j}+\mu=(B\odot A)_{i,j}$.

For the matrix 
$$
C=\left[\begin{matrix}\begin{matrix} \mu & \varepsilon\\ \varepsilon & \mu \end{matrix} & - \infty \\
  -\infty & \mu I_{n-2} \end{matrix}\right], 
$$ 
it is easy to see that $B\odot C=C\odot B$. Therefore, 
$$
D-  \left[\begin{matrix}d_{1,1}  I_2 & - \infty \\
  -\infty & d_{2,2}I_{n-2} \end{matrix}\right] - C- B-A$$ is a path in $\Gamma(\mm_n(\TT))$ and thus $d(D,A) \leq 4$.
%
% Guterman:
%Let $n\ge 3$. We denote by $J_n$ the matrix of all ones (in $R_\max$, i.e. $0$ in $R$) and proof that the distance from $J_n$ 
%to arbitrary matrix is $\le 3$. We split the proof into 3 cases.
%1. $A$ is not diagonal.
%Let $\lambda=\max_{i,j} a_{i,j}$. Then $[A,A \oplus \lambda \odot I_n]=0$. 
%As far as $(A \oplus \lambda \odot I_n)\odot J_n=J_n \odot (A \oplus \lambda \odot I_n)$ is the matrix of all
%$\lambda$.
%Thus $A$ --- $A \oplus \lambda \odot I_n$ --- $J_n$ is a path.
%2. $A$ is diagonal with some nonzero $\lambda_1,\ldots,\lambda_k$ on the diagonal. WLOG (up to a permutation)
%$A=\lambda_1 \odot E_{11} \oplus \ldots \oplus \lambda_k \odot E_{kk}$ for some $k\le n-2$. 
%Then $A$ commutes with the block matrix with the blocks $I_k$ and $J_{n-k}$. And the last one commutes with $J_n$.
%3. Up to a permutation $A=\lambda_1 \odot E_{11} \oplus \ldots \oplus \lambda_{n-1} \odot E_{n-1n-1}$.
%Let $\lambda=\max_i \lambda_i$. Then $A$ --- $A \oplus \lambda \odot I_{n-1}$  --- block matrix with blocks
%$J_{n-1}$ and $I_1$ --- $J_n$ is a path. 
\end{proof}

\bigskip
As a Corollary of \cite[Thm.~7]{DolObl10} we have the following.

\medskip

\begin{Proposition}\label{thm:glntropical}
For the tropical semiring $\TT$, we have 
\begin{equation*}
\diam{\Gamma({\rm GL}_{n}(\TT))}=\begin{cases}
 5, & n \text{ prime},\\
 \infty, & \text{otherwise.}
 \end{cases} 
 \end{equation*}
\end{Proposition}

\bigskip

\begin{Remark}
Note that  \cite[Thm.~14, 15]{AMRR06} state that for an infinite commutative division ring $R$, we have the 
equality
$$\diam{\Gamma(\mm_n(R))}=\diam{\Gamma({\rm GL}_n(R))}$$
for all $n \geq 2$.
With  Theorem \ref{thm:mntropical} and Proposition \ref{thm:glntropical}, we proved that this result cannot
be generalized to commutative division semirings. 
\end{Remark}

\bigskip
\bigskip
 
%-----------------------------------------------------
%-----------------------------------------------------
\section{The commuting graph of the full matrix semiring over a nonentire commutative semiring}
%-----------------------------------------------------
%-----------------------------------------------------

\bigskip

 In this section, we prove that the diameter of the commuting graph of the full matrix semiring over a nonentire commutative semiring
 is always equal to 3.  This generalizes the result from \cite[Thm.~1.1]{giudici}, where it has been proved that the diameter
 of the commuting graph of the full matrix ring over the ring ${\mathbb Z}_m$ for a composite number $m$ is always equal to 3.

 \bigskip

First, we find the lower bound for the diameter of the commuting graph of matrices over an
arbitrary commutative semiring.

 \bigskip
 
 \begin{Proposition}\label{thm:3}
   If $S$ is a commutative semiring, then $\diam{\Gamma(\mm_n(S))}\geq 3$
   for $n\geq 2$. 
 \end{Proposition}

 \medskip

 \begin{proof}
   By Lemma \ref{lem:centrJ}, the only matrices in the intersection of the centralizer of $J_n$ and
   the centralizer of $J_n^T$ are scalar matrices, but these are of course central.
   Thus, the distance between matrices $J_n$ and $J_n^T$ in $\Gamma(\mm_n(S))$ is at least $3$.
 \end{proof}

 \bigskip 

We now prove that this bound is sharp, since the diameter is equal to 
3 for all commuting graphs of matrices over nonentire commutative semirings. 
 
 \bigskip 
 
 \begin{Theorem}
   If $S$ is a nonentire commutative semiring, then $\diam{\Gamma(\mm_n(S))} = 3$
   for $n \geq 2$.
 \end{Theorem}
 
 \medskip
 
 \begin{proof}
   Since $S$ is not entire, there exist nonzero $ x, y \in S$ such that $xy=yx=0$. (Note that it may happen that $x=y$.)
   Choose arbitrary noncentral matrices $A, B \in \mm_n(S)$.  
   \begin{enumerate}[{Case} 1:]
   \item
   Suppose $xA$ and $yB$ are noncentral matrices in $\mm_n(S)$. 
   Then $A-xA-yB-B$ is a path of length $3$ in $\Gamma(\mm_n(S))$.  
   \item
   Suppose $xA$ is central and $yB$ is a noncentral matrix in $\mm_n(S)$. 
   Then $AxE_{1,2}=xAE_{1,2}=E_{1,2}xA=xE_{1,2}A$, so
   $A-xE_{1,2}-yB-B$ is a path of length $3$ in $\Gamma(\mm_n(S))$.  
   Similarly, we treat the case when $yB$ is central and $xA$ is not.
   \item
   Suppose $xA$ and $yB$ are central matrices in $\mm_n(S)$. 
   Then $A-xE_{1,2}-yE_{1,2}-B$ is a path of length $3$ in $\Gamma(\mm_n(S))$.  
   \end{enumerate}
   The three cases considered show that the distance between any two matrices in $\mm_n(S)$
   is at most $3$, and therefore $\diam{\Gamma(\mm_n(S))}\leq 3$.
   By Proposition \ref{thm:3}, it follows that $\diam{\Gamma(\mm_n(S))} = 3$.
 \end{proof}

\end{document}